\documentclass{amsart}
\usepackage{amsmath}
\usepackage{amsfonts}

\newtheorem{theorem}{Theorem}
\theoremstyle{plain}

\newtheorem{corollary}{Corollary}

\numberwithin{equation}{section}
\input{tcilatex}

\begin{document}
\title[$q$-Genocchi numbers and polynomials with weight $\alpha $ and $\beta 
$]{A note on the values of the weighted $q$-Bernstein Polynomials and
modified$\ q$-Genocchi Numbers with weight $\alpha $ and $\beta $ via the $p$%
-adic $q$-integral on $%
\mathbb{Z}
_{p}$}
\author{Serkan Arac\i }
\address{University of Gaziantep, Faculty of Science and Arts, Department of
Mathematics, 27310 Gaziantep, TURKEY}
\email{mtsrkn@hotmail.com}
\author{Mehmet Acikgoz}
\address{University of Gaziantep, Faculty of Science and Arts, Department of
Mathematics, 27310 Gaziantep, TURKEY}
\email{acikgoz@gantep.edu.tr}
\date{January 8, 2012}

\begin{abstract}
The rapid development of $q$-calculus has led to the discovery of new
generalizations of Bernstein polynomials and Genocchi polynomials involving $%
q$-integers. The present paper deals with weighted $q$-Bernstein polynomials
and $q$-Genocchi numbers with weight $\alpha $ and $\beta $. We apply the
method of generating function and $p$-adic $q$-integral representation on $%
\mathbb{Z}
_{p}$, which are exploited to derive further classes of Bernstein
polynomials and $q$-Genocchi numbers and polynomials. To be more precise we
summarize our results as follows, we obtain some combinatorial relations
between $q$-Genocchi numbers and polynomials with weight $\alpha $ and $%
\beta $. Furthermore we derive an integral representation of weighted $q$%
-Bernstein polynomials of degree $n$ on $%
\mathbb{Z}
_{p}$. Also we deduce a fermionic $p$-adic $q$-integral representation of
product weighted $q$-Bernstein polynomials of different degrees $%
n_{1},n_{2},\cdots $\ on $%
\mathbb{Z}
_{p}$ and show that it can be written with $q$-Genocchi numbers with weight $%
\alpha $ and $\beta $ which yields a deeper insight into the effectiveness
of this type of generalizations. Our new generating function possess a
number of interesting properties which we state in this paper
\end{abstract}

\maketitle

\section{Introduction, Definitions and Notations}

The $q$-calculus theory is a novel theory that is based on finite difference
re-scaling. First results in $q$-calculus belong to Euler, who discovered
Euler's Identities for $q$-exponential functions and Gauss, who discovered $%
q $-binomial formula. The systematic development of $q$-calculus begins from
F. H. Jackson who 1908 reintroduced the Euler Jackson $q$-difference
operator (Jackson, 1908). One of important branches of $q$-calculus is $q$%
-type special orthogonal polynomials. Also $p$-adic numbers were invented by
Kurt Hensel around the end of the nineteenth century and these two branches
of number theory jointed with the link of $p$-adic $q$-integral and
developed. In spite of their being already one hundred years old, these
special numbers and polynomials, for instance $q$-Bernstein numbers and
polynomials, $q$-Genocchi numbers and polynomials and etc. are still today
enveloped in an aura of mystery within the scientific community. The $p$%
-adic integral was used in mathematical physics, for instance, the
functional equation of the $q$-zeta function, $q$-stirling numbers and $q$%
-Mahler theory of integration with respect to the ring $%
\mathbb{Z}
_{p}$ together with Iwasawa's $p$-adic $q$-$L$ functions. Professor T. Kim 
\cite{Kim 21}, also studied on $p$-adic interpolation functions of special
orthogonal polynomials. In during the last ten years, the $q$-Bernstein
polynomials and $q$-Genocchi polynomials have attracted a lot of interest
and have been studied from different angles along with some generalizations
and modifications by a number of researchers. By using the $p$-adic
invariant $q$-integral on $%
\mathbb{Z}
_{p}$, Professor T. Kim in \cite{kim 18}, constructed $p$-adic Bernoulli
numbers and polynomials with weight $\alpha $. After Seo and first author in 
\cite{Seo}, extended Kim's method for $q$-Genocchi numbers and polynomials
and also they defined $q$-Genocchi numbers and polynomials with weight $%
\alpha $ and $\beta $. Our aim of this paper is to show that a fermionic $p$%
-adic $q$-integral representation of product weighted $q$-Bernstein
polynomials of different degrees $n_{1},n_{2},\cdots $ on $%
\mathbb{Z}
_{p}$ can be written with $q$-Genocchi numbers with weight $\alpha $ and $%
\beta $.

Let $p$ be a fixed odd prime number. Throughout this paper we use the
following notations. By $%
\mathbb{Z}
_{p}$ we denote the ring of $p$-adic rational integers, $%
\mathbb{Q}
$ denotes the field of rational numbers, $%
\mathbb{Q}
_{p}$ denotes the field of $p$-adic rational numbers, and $%
\mathbb{C}
_{p}$ denotes the completion of algebraic closure of $%
\mathbb{Q}
_{p}$. Let $%
\mathbb{N}
$ be the set of natural numbers and $%
\mathbb{N}
^{\ast }=%
\mathbb{N}
\cup \left\{ 0\right\} $. The $p$-adic absolute value is defined by $%
\left\vert p\right\vert _{p}=\frac{1}{p}$. In this paper we assume $%
\left\vert q-1\right\vert _{p}<1$ as an indeterminate. In [23-25], let $%
UD\left( 
\mathbb{Z}
_{p}\right) $ be the space of uniformly differentiable functions on $%
\mathbb{Z}
_{p}$. For $f\in UD\left( 
\mathbb{Z}
_{p}\right) $, the fermionic $p$-adic $q$-integral on $%
\mathbb{Z}
_{p}$ is defined by T. Kim:%
\begin{eqnarray}
I_{-q}\left( f\right) &=&\int_{%
\mathbb{Z}
_{p}}f\left( \xi \right) d\mu _{-q}\left( \xi \right)  \label{equation 1} \\
&=&\lim_{N\rightarrow \infty }\frac{1}{\left[ p^{N}\right] _{-q}}\sum_{\xi
=0}^{p^{N}-1}q^{\xi }f\left( \xi \right) \left( -1\right) ^{\xi }.\text{ } 
\notag
\end{eqnarray}

For $\alpha ,k,n\in 
\mathbb{N}
^{\ast }$ and $x\in \left[ 0,1\right] $, T. Kim et al$.$ defined weighted $q$%
-Bernstein polynomials as follows:%
\begin{equation}
B_{k,n}^{\left( \alpha \right) }\left( x,q\right) =\binom{n}{k}\left[ x%
\right] _{q^{\alpha }}^{k}\left[ 1-x\right] _{q^{-\alpha }}^{n-k}\text{, \
(for detail, see [3, 27, 33, 34]). }  \label{equation 2}
\end{equation}

In (\ref{equation 2}), we put $q\rightarrow 1$ and $\alpha =1$, $\left[ x%
\right] _{q^{\alpha }}^{k}\rightarrow x^{k},$ $\left[ 1-x\right]
_{q^{-\alpha }}^{n-k}\rightarrow \left( 1-x\right) ^{n-k}$ and we obtain the
classical Bernstein polynomials (see\cite{acikgoz2}, \cite{acikgoz 4}),

where, $\left[ x\right] _{q}$ is a $q$-extension of $x$ which is defined by%
\begin{equation*}
\left[ x\right] _{q}=\frac{1-q^{x}}{1-q}\text{, \ \ (see [2-28, 32-34]).}
\end{equation*}

Note that $\lim_{q\rightarrow 1}\left[ x\right] _{q}=x$.

In previous paper [8], for $n\in 
\mathbb{N}
^{\ast }$, modified $q$-Genocchi numbers with weight $\alpha $ and $\beta $
are defined by Araci et al. as follows:%
\begin{eqnarray}
\frac{g_{n+1,q}^{\left( \alpha ,\beta \right) }\left( x\right) }{n+1}
&=&\int_{%
\mathbb{Z}
_{p}}q^{-\beta \xi }\left[ x+\xi \right] _{q^{\alpha }}^{n}d\mu _{-q^{\beta
}}\left( \xi \right)  \label{equation 3} \\
&=&\frac{\left[ 2\right] _{q^{\beta }}}{\left[ \alpha \right] _{q}^{n}\left(
1-q\right) ^{n}}\sum_{l=0}^{n}\binom{n}{l}\left( -1\right) ^{l}q^{\alpha
\ell x}\frac{1}{1+q^{\alpha \ell }}  \notag \\
&=&\left[ 2\right] _{q^{\beta }}\sum_{m=0}^{\infty }\left( -1\right) ^{m}%
\left[ m+x\right] _{q^{\alpha }}^{n}.  \notag
\end{eqnarray}

In the special case, $x=0$, $g_{n,q}^{\left( \alpha ,\beta \right) }\left(
0\right) =g_{n,q}^{\left( \alpha ,\beta \right) }$ are called the $q$%
-Genocchi numbers with weight $\alpha $ and $\beta $.

In [8], for $\alpha \in 
\mathbb{N}
^{\ast }$ and $n\in 
\mathbb{N}
$, $q$-Genocchi numbers with weight $\alpha $ and $\beta $ are defined by
Araci et al. as follows:%
\begin{equation}
g_{0,q}^{\left( \alpha ,\beta \right) }=0,\text{ and }g_{n,q}^{\left( \alpha
,\beta \right) }\left( 1\right) +g_{n,q}^{\left( \alpha ,\beta \right)
}=\left\{ \QATOP{\left[ 2\right] _{q^{\beta }},\text{ if }n=1,}{0,\text{ \ \
\ \ if }n>1.}\right.  \label{equation 18}
\end{equation}

\qquad In this paper, we obtained some relations between the weighted $q$%
-Bernstein polynomials and the modified $q$-Genocchi numbers with weight $%
\alpha $ and $\beta $. From these relations, we derive some interesting
identities on the $q$-Genocchi numbers with weight $\alpha $ and $\beta $.

\section{On the $q$-Genocchi numbers and polynomials with weight $\protect%
\alpha $ and $\protect\beta $}

By the definition of $q$-Genocchi polynomials with weight $\alpha $ and $%
\beta $, we easily get 
\begin{eqnarray*}
\frac{g_{n+1,q}^{\left( \alpha ,\beta \right) }\left( x\right) }{n+1}
&=&\int_{%
\mathbb{Z}
_{p}}q^{-\beta \xi }\left[ x+\xi \right] _{q^{\alpha }}^{n}d\mu _{-q^{\beta
}}\left( \xi \right) \\
&=&\int_{%
\mathbb{Z}
_{p}}q^{-\beta \xi }\left( \left[ x\right] _{q^{\alpha }}+q^{\alpha x}\left[
\xi \right] _{q^{\alpha }}\right) ^{n}d\mu _{-q}\left( \xi \right) \\
&=&\sum_{k=0}^{n}\binom{n}{k}\left[ x\right] _{q^{\alpha }}^{n-k}q^{\alpha
kx}\int_{%
\mathbb{Z}
_{p}}q^{-\beta \xi }\left[ \xi \right] _{q^{\alpha }}^{k}d\mu _{-q}\left(
\xi \right) \\
&=&\sum_{k=0}^{n}\binom{n}{k}\left[ x\right] _{q^{\alpha }}^{n-k}q^{\alpha
kx}\frac{g_{k+1,q}^{\left( \alpha ,\beta \right) }}{k+1}.
\end{eqnarray*}

Therefore, we obtain the following Theorem:

\begin{theorem}
For $n,\alpha ,\beta \in 
\mathbb{N}
^{\ast },$ we have%
\begin{equation}
g_{n,q}^{\left( \alpha ,\beta \right) }\left( x\right) =q^{-\alpha
x}\sum_{k=0}^{n}\binom{n}{k}q^{\alpha kx}g_{k,q}^{\left( \alpha ,\beta
\right) }\left[ x\right] _{q^{\alpha }}^{n-k},  \label{equation 6}
\end{equation}%
Moreover,%
\begin{equation}
g_{n,q}^{\left( \alpha ,\beta \right) }\left( x\right) =q^{-\alpha x}\left(
q^{\alpha x}g_{q}^{\left( \alpha ,\beta \right) }+\left[ x\right]
_{q^{\alpha }}\right) ^{n},  \label{equation 13}
\end{equation}%
by using the $umbral$(symbolic) convention $\left( g_{q}^{\left( \alpha
,\beta \right) }\right) ^{n}=g_{n,q}^{\left( \alpha ,\beta \right) }.$
\end{theorem}

By expression of (\ref{equation 3}), we get 
\begin{eqnarray*}
\frac{g_{n+1,q^{-1}}^{\left( \alpha ,\beta \right) }\left( 1-x\right) }{n+1}
&=&\int_{%
\mathbb{Z}
_{p}}q^{\beta \xi }\left[ 1-x+\xi \right] _{q^{-\alpha }}^{n}d\mu
_{-q^{-\beta }}\left( \xi \right) \\
&=&\frac{\left[ 2\right] _{q^{-\beta }}}{\left( 1-q^{-\alpha }\right) ^{n}}%
\sum_{l=0}^{n}\binom{n}{l}\left( -1\right) ^{l}q^{-\alpha \ell \left(
1-x\right) }\frac{1}{1+q^{-\alpha \ell }} \\
&=&\left( -1\right) ^{n}q^{\alpha n-\beta }\left( \frac{\left[ 2\right]
_{q^{\beta }}}{\left( 1-q^{\alpha }\right) ^{n}}\sum_{l=0}^{n}\binom{n}{l}%
\left( -1\right) ^{l}q^{\alpha lx}\frac{1}{1+q^{\alpha l}}\right) \\
&=&\left( -1\right) ^{n}q^{\alpha n-\beta }\frac{g_{n+1,q}^{\left( \alpha
,\beta \right) }\left( x\right) }{n+1}.
\end{eqnarray*}

Consequently, we obtain the following Theorem:

\begin{theorem}
The following 
\begin{equation}
g_{n+1,q^{-1}}^{\left( \alpha ,\beta \right) }\left( 1-x\right) =\left(
-1\right) ^{n}q^{\alpha n-\beta }g_{n+1,q}^{\left( \alpha ,\beta \right)
}\left( x\right)  \label{equation 5}
\end{equation}%
is true.
\end{theorem}

From expression of (\ref{equation 13}) and Theorem 1, we get the following
Theorem:

\begin{theorem}
The following identity holds 
\begin{equation*}
g_{0,q}^{\left( \alpha ,\beta \right) }=0,\text{ and }q^{-\alpha }\left(
q^{\alpha }g_{q}^{\left( \alpha ,\beta \right) }+1\right)
^{n}+g_{n,q}^{\left( \alpha ,\beta \right) }=\left\{ \QATOP{\left[ 2\right]
_{q^{\beta }},\text{ if }n=1,}{0,\text{ \ \ \ \ \ if }n>1,}\right.
\end{equation*}%
with the usual convention about replacing $\left( g_{q}^{\left( \alpha
,\beta \right) }\right) ^{n}$ by $g_{n,q}^{\left( \alpha ,\beta \right) }$.
\end{theorem}

\qquad For $n,\alpha \in 
\mathbb{N}
$, by Theorem 3, we note that 
\begin{eqnarray*}
q^{2\alpha }g_{n,q}^{\left( \alpha ,\beta \right) }\left( 2\right) &=&\left(
q^{\alpha }\left( q^{\alpha }g_{q}^{\left( \alpha ,\beta \right) }+1\right)
+1\right) ^{n} \\
&=&\sum_{k=0}^{n}\binom{n}{k}q^{k\alpha }\left( q^{\alpha }g_{q}^{\left(
\alpha ,\beta \right) }+1\right) ^{k} \\
&=&\left( q^{\alpha }g_{q}^{\left( \alpha ,\beta \right) }+1\right)
^{0}+nq^{\alpha }\left( q^{\alpha }g_{q}^{\left( \alpha ,\beta \right)
}+1\right) ^{1} \\
&&+\sum_{k=2}^{n}\binom{n}{k}q^{k\alpha }\left( q^{\alpha }g_{q}^{\left(
\alpha ,\beta \right) }+1\right) ^{k} \\
&=& n q^{2\alpha }\left[ 2\right] _{q^{\beta }}-q^{\alpha }\sum_{k=0}^{n}%
\binom{n}{k}q^{\alpha k}g_{k,q}^{\left( \alpha ,\beta \right) } \\
&=& n q^{2 \alpha }\left[ 2\right] _{q^{\beta }}+q^{\alpha }g_{n,q}^{\left(
\alpha ,\beta \right) },\text{ if }n>1.
\end{eqnarray*}

Consequently, we state the following Theorem:

\begin{theorem}
For $n\in 
\mathbb{N}
$, we have 
\begin{equation*}
g_{n,q}^{\left( \alpha ,\beta \right) }\left( 2\right) ={\ n \left[ 2\right]
_{q^{\beta }}}+ \dfrac{g_{n,q}^{\left( \alpha ,\beta \right) }}{q^{\alpha }}.
\end{equation*}
\end{theorem}

From expression of\ Theorem 2 and (\ref{equation 5}), we easily see that%
\begin{eqnarray}
&&\left( n+1\right) q^{-\beta }\int_{%
\mathbb{Z}
_{p}}q^{-\beta \xi }\left[ 1-\xi \right] _{q^{-\alpha }}^{n}d\mu _{-q^{\beta
}}\left( \xi \right)  \label{equation 7} \\
&=&\left( -1\right) ^{n}q^{n\alpha -\beta }\int_{%
\mathbb{Z}
_{p}}q^{-\beta \xi }\left[ \xi -1\right] _{q^{\alpha }}^{n}d\mu _{-q^{\beta
}}\left( \xi \right)  \notag \\
&=&\left( -1\right) ^{n}q^{n\alpha -\beta }g_{n+1,q}^{\left( \alpha ,\beta
\right) }\left( -1\right) =g_{n+1,q^{-1}}^{\left( \alpha ,\beta \right)
}\left( 2\right) .  \notag
\end{eqnarray}

Thus, we obtain the following Theorem.

\begin{theorem}
The following identity 
\begin{equation*}
\left( n+1\right) q^{-\beta }\int_{%
\mathbb{Z}
_{p}}q^{-\beta \xi }\left[ 1-\xi \right] _{q^{-\alpha }}^{n}d\mu _{-q^{\beta
}}\left( \xi \right) =g_{n+1,q^{-1}}^{\left( \alpha ,\beta \right) }\left(
2\right)
\end{equation*}%
is true.
\end{theorem}

Let $n,\alpha \in 
\mathbb{N}
$. By expression of Theorem 4 and Theorem 5, we get%
\begin{eqnarray}
&&\left( n+1\right) q^{-\beta }\int_{%
\mathbb{Z}
_{p}}q^{-\beta \xi }\left[ 1-\xi \right] _{q^{-\alpha }}^{n}d\mu _{-q^{\beta
}}\left( \xi \right)  \label{equation 8} \\
&=&\left( {n+1}\right) q^{-\beta }{\left[ 2\right] _{q^{\beta }}}+q^{\alpha
}g_{n+1,q^{-1}}^{\left( \alpha ,\beta \right) }.  \notag
\end{eqnarray}

For (\ref{equation 8}), we obtain corollary as follows:

\begin{corollary}
For $n,\alpha \in 
\mathbb{N}
^{\ast },$ we have%
\begin{equation*}
\int_{%
\mathbb{Z}
_{p}}q^{-\beta \xi }\left[ 1-\xi \right] _{q^{-\alpha }}^{n}d\mu _{-q^{\beta
}}\left( \xi \right) =\left[ 2\right] _{q^{\beta }}+q^{\alpha -\beta }\frac{%
g_{n+1,q^{-1}}^{\left( \alpha ,\beta \right) }}{n+1}.
\end{equation*}
\end{corollary}

\section{Novel identities on the weighted $q$-Genocchi numbers}

\qquad In this section, we develop modifed $q$-Genocchi numbers with weight $%
\alpha $ and $\beta $, namely, we derive interesting and worthwhile
relations in Analytic Number Theory.

\qquad For $x\in 
\mathbb{Z}
_{p}$, the $p$-adic analogues of weighted $q$-Bernstein polynomials are
given by%
\begin{equation}
B_{k,n}^{\left( \alpha \right) }\left( x,q\right) =\binom{n}{k}\left[ x%
\right] _{q^{\alpha }}^{k}\left[ 1-x\right] _{q^{-\alpha }}^{n-k},\text{
where }n,k,\alpha \in 
\mathbb{N}
^{\ast }.  \label{equation 9}
\end{equation}

By expression of (\ref{equation 9}), Kim et. al. get the symmetry of $q$%
-Bernstein polynomials weight $\alpha $ as follows:%
\begin{equation}
B_{k,n}^{\left( \alpha \right) }\left( x,q\right) =B_{n-k,n}^{\left( \alpha
\right) }\left( 1-x,q^{-1}\right) \text{, (for detail, see \cite{kim 19}).}
\label{equation 10}
\end{equation}

Thus, from Corollary 1, (\ref{equation 9}) and (\ref{equation 10}), we see
that%
\begin{eqnarray*}
\int_{%
\mathbb{Z}
_{p}}B_{k,n}^{\left( \alpha \right) }\left( \xi ,q\right) q^{-\beta \xi
}d\mu _{-q^{\beta }}\left( \xi \right) &=&\int_{%
\mathbb{Z}
_{p}}B_{n-k,n}^{\left( \alpha \right) }\left( 1-\xi ,q^{-1}\right) q^{-\beta
\xi }d\mu _{-q^{\beta }}\left( \xi \right) \\
&=&\binom{n}{k}\sum_{l=0}^{k}\binom{k}{l}\left( -1\right) ^{k+l}\int_{%
\mathbb{Z}
_{p}}q^{-\beta \xi }\left[ 1-\xi \right] _{q^{-\alpha }}^{n-l}d\mu
_{-q^{\beta }}\left( \xi \right) \\
&=&\binom{n}{k}\sum_{l=0}^{k}\binom{k}{l}\left( -1\right) ^{k+l}\left( \left[
2\right] _{q^{\beta }}+q^{\alpha -\beta }\frac{g_{n-l+1,q^{-1}}^{\left(
\alpha ,\beta \right) }}{n-l+1}\right) .
\end{eqnarray*}

For $n$, $k\in 
\mathbb{N}
^{\ast }$ and $\alpha \in 
\mathbb{N}
$ with $n>k$, we obtain%
\begin{eqnarray}
&&\int_{%
\mathbb{Z}
_{p}}B_{k,n}^{\left( \alpha \right) }\left( \xi ,q\right) q^{-\beta \xi
}d\mu _{-q^{\beta }}\left( \xi \right)  \notag \\
&=&\binom{n}{k}\sum_{l=0}^{k}\binom{k}{l}\left( -1\right) ^{k+l}\left( \left[
2\right] _{q^{\beta }}+q^{\alpha -\beta }\frac{g_{n-l+1,q^{-1}}^{\left(
\alpha ,\beta \right) }}{n-l+1}\right)  \label{equation 11} \\
&=&\left\{ \QATOPD. . {\left[ 2\right] _{q^{\beta }}+q^{\alpha -\beta }\frac{%
g_{n+1,q^{-1}}^{\left( \alpha ,\beta \right) }}{n+1},\text{ \ \ \ \ \ \ \ \
\ \ \ \ \ \ \ \ \ \ \ \ \ \ \ \ \ \ \ \ \ \ \ \ \ \ \ if }k=0,}{\binom{n}{k}%
\sum_{l=0}^{k}\binom{k}{l}\left( -1\right) ^{k+l}\left( \left[ 2\right]
_{q^{\beta }}+q^{\alpha -\beta }\frac{g_{n-l+1,q^{-1}}^{\left( \alpha ,\beta
\right) }}{n-l+1}\right) ,\text{ if }k>0.}\right.  \notag
\end{eqnarray}

Let us take the fermionic $p$-adic $q$-integral on $%
\mathbb{Z}
_{p}$ on the weighted $q$-Bernstein polynomials of degree $n$ as follows:%
\begin{eqnarray}
\int_{%
\mathbb{Z}
_{p}}B_{k,n}^{\left( \alpha \right) }\left( \xi ,q\right) q^{-\beta \xi
}d\mu _{-q^{\beta }}\left( \xi \right) &=&\binom{n}{k}\int_{%
\mathbb{Z}
_{p}}q^{-\beta \xi }\left[ \xi \right] _{q^{\alpha }}^{k}\left[ 1-\xi \right]
_{q^{-\alpha }}^{n-k}d\mu _{-q^{\beta }}\left( \xi \right)
\label{equation 12} \\
&=&\binom{n}{k}\sum_{l=0}^{n-k}\binom{n-k}{l}\left( -1\right) ^{l}\frac{%
g_{l+k+1,q}^{\left( \alpha ,\beta \right) }}{l+k+1}.  \notag
\end{eqnarray}

Consequently, by expression of (\ref{equation 11}) and (\ref{equation 12})$,$
we state the following Theorem:

\begin{theorem}
The following identity holds%
\begin{equation*}
\sum_{l=0}^{n-k}\binom{n-k}{l}\left( -1\right) ^{l}\frac{g_{l+k+1,q}^{\left(
\alpha ,\beta \right) }}{l+k+1}=\left\{ \QATOPD. . {\left[ 2\right]
_{q^{\beta }}+q^{\alpha -\beta }\frac{g_{n+1,q^{-1}}^{\left( \alpha ,\beta
\right) }}{n+1},\text{ \ \ \ \ \ \ \ \ \ \ \ \ \ \ \ \ \ \ \ \ \ \ \ \ \ \ \
\ if }k=0,}{\sum_{l=0}^{k}\binom{k}{l}\left( -1\right) ^{k+l}\left( \left[ 2%
\right] _{q^{\beta }}+q^{\alpha -\beta }\frac{g_{n-l+1,q^{-1}}^{\left(
\alpha ,\beta \right) }}{n-l+1}\right) ,\text{ if }k>0.}\right.
\end{equation*}
\end{theorem}

Let $n_{1},n_{2},k\in 
\mathbb{N}
^{\ast }$ and $\alpha \in 
\mathbb{N}
$ with $n_{1}+n_{2}>2k$. Then, we get%
\begin{eqnarray}
&&\int_{%
\mathbb{Z}
_{p}}B_{k,n_{1}}^{\left( \alpha \right) }\left( \xi ,q\right)
B_{k,n_{2}}^{\left( \alpha \right) }\left( \xi ,q\right) q^{-\beta \xi }d\mu
_{-q^{\beta }}\left( \xi \right)  \notag \\
&=&\binom{n_{1}}{k}\binom{n_{2}}{k}\sum_{l=0}^{2k}\binom{2k}{l}\left(
-1\right) ^{2k+l}\int_{%
\mathbb{Z}
_{p}}q^{-\beta \xi }\left[ 1-\xi \right] _{q^{-\alpha }}^{n_{1}+n_{2}-l}d\mu
_{-q^{\beta }}\left( \xi \right)  \notag \\
&=&\left( \binom{n_{1}}{k}\binom{n_{2}}{k}\sum_{l=0}^{2k}\binom{2k}{l}\left(
-1\right) ^{2k+l}\left( \left[ 2\right] _{q^{\beta }}+q^{\alpha -\beta }%
\frac{g_{n_{1}+n_{2}-l+1,q^{-1}}^{\left( \alpha ,\beta \right) }}{%
n_{1}+n_{2}-l+1}\right) \right)  \notag \\
&=&\left\{ \QATOPD. . {\left[ 2\right] _{q^{\beta }}+q^{\alpha -\beta }\frac{%
g_{n_{1}+n_{2}+1,q^{-1}}^{\left( \alpha ,\beta \right) }}{n_{1}+n_{2}+1},%
\text{ \ \ \ \ \ \ \ \ \ \ \ \ \ \ \ \ \ \ \ \ \ \ \ \ \ \ \ \ \ \ \ \ \ \ \
\ \ \ \ \ \ \ \ \ \ \ \ if }k=0,}{\binom{n_{1}}{k}\binom{n_{2}}{k}%
\sum_{l=0}^{2k}\binom{2k}{l}\left( -1\right) ^{2k+l}\left( \left[ 2\right]
_{q^{\beta }}+q^{\alpha -\beta }\frac{g_{n_{1}+n_{2}-l+1,q^{-1}}^{\left(
\alpha ,\beta \right) }}{n_{1}+n_{2}-l+1}\right) ,\text{ \ if }k\neq
0.}\right.  \notag
\end{eqnarray}

Therefore, we obtain the following Theorem:

\begin{theorem}
For $n_{1},n_{2},k\in 
\mathbb{N}
^{\ast }$and $\alpha ,\beta \in 
\mathbb{N}
$ with $n_{1}+n_{2}>2k,$ we have%
\begin{eqnarray*}
&&\int_{%
\mathbb{Z}
_{p}}q^{-\beta \xi }B_{k,n_{1}}^{\left( \alpha \right) }\left( \xi ,q\right)
B_{k,n_{2}}^{\left( \alpha \right) }\left( \xi ,q\right) d\mu _{-q^{\beta
}}\left( \xi \right) \\
&=&\left\{ \QATOPD. . {\left[ 2\right] _{q^{\beta }}+q^{\alpha -\beta }\frac{%
g_{n_{1}+n_{2}+1,q^{-1}}^{\left( \alpha ,\beta \right) }}{n_{1}+n_{2}+1},%
\text{ \ \ \ \ \ \ \ \ \ \ \ \ \ \ \ \ \ \ \ \ \ \ \ \ \ \ \ \ \ \ \ \ \ \ \
\ \ if }k=0,}{\binom{n_{1}}{k}\binom{n_{2}}{k}\sum_{l=0}^{2k}\binom{2k}{l}%
\left( -1\right) ^{2k+l}\left( \left[ 2\right] _{q^{\beta }}+q^{\alpha
-\beta }\frac{g_{n_{1}+n_{2}-l+1,q^{-1}}^{\left( \alpha ,\beta \right) }}{%
n_{1}+n_{2}-l+1}\right) ,\text{ \ if }k\neq 0.}\right.
\end{eqnarray*}
\end{theorem}

By using the binomial theorem, we can derive the following equation.%
\begin{eqnarray}
&&\int_{%
\mathbb{Z}
_{p}}B_{k,n_{1}}^{\left( \alpha \right) }\left( \xi ,q\right)
B_{k,n_{2}}^{\left( \alpha \right) }\left( \xi ,q\right) q^{-\beta \xi }d\mu
_{-q^{\beta }}\left( \xi \right)  \notag \\
&=&\dprod\limits_{i=1}^{2}\binom{n_{i}}{k}\sum_{l=0}^{n_{1}+n_{2}-2k}\binom{%
n_{1}+n_{2}-2k}{l}\left( -1\right) ^{l}\int_{%
\mathbb{Z}
_{p}}\left[ \xi \right] _{q^{\alpha }}^{2k+l}q^{-\beta \xi }d\mu _{-q^{\beta
}}\left( \xi \right)  \label{equation 14} \\
&=&\dprod\limits_{i=1}^{2}\binom{n_{i}}{k}\sum_{l=0}^{n_{1}+n_{2}-2k}\binom{%
n_{1}+n_{2}-2k}{l}\left( -1\right) ^{l}\frac{g_{l+2k+1,q}^{\left( \alpha
,\beta \right) }}{l+2k+1}.  \notag
\end{eqnarray}

Thus, we can obtain the following Corollary:

\begin{corollary}
For $n_{1},n_{2},k\in 
\mathbb{N}
^{\ast }$ and $\alpha \in 
\mathbb{N}
$ with $n_{1}+n_{2}>2k,$ we have%
\begin{eqnarray*}
&&\sum_{l=0}^{n_{1}+n_{2}-2k}\binom{n_{1}+n_{2}-2k}{l}\left( -1\right) ^{l}%
\frac{g_{l+2k+1,q}^{\left( \alpha ,\beta \right) }}{l+2k+1} \\
&=&\left\{ \QATOPD. . {\left[ 2\right] _{q^{\beta }}+q^{\alpha -\beta }\frac{%
g_{n_{1}+n_{2}+1,q^{-1}}^{\left( \alpha ,\beta \right) }}{n_{1}+n_{2}+1},%
\text{ \ \ \ \ \ \ \ \ \ \ \ \ \ \ \ \ \ \ \ \ \ \ \ \ \ \ \ \ if }%
k=0,}{\sum_{l=0}^{2k}\binom{2k}{l}\left( -1\right) ^{2k+l}\left( \left[ 2%
\right] _{q^{\beta }}+q^{\alpha -\beta }\frac{g_{n_{1}+n_{2}-l+1,q^{-1}}^{%
\left( \alpha ,\beta \right) }}{n_{1}+n_{2}-l+1}\right) ,\text{ \ if }k\neq
0.}\right.
\end{eqnarray*}
\end{corollary}

For $\xi \in 
\mathbb{Z}
_{p}$ and $s\in 
\mathbb{N}
$ with $s\geq 2,$ let $n_{1},n_{2},...,n_{s},k\in 
\mathbb{N}
^{\ast }$ and $\alpha \in 
\mathbb{N}
$ with $\sum_{l=1}^{s}n_{l}>sk$. Then we take the fermionic $p$-adic $q$%
-integral on $%
\mathbb{Z}
_{p}$ for the weighted $q$-Bernstein polynomials of degree $n$ as follows:%
\begin{eqnarray*}
&&\int_{%
\mathbb{Z}
_{p}}\underset{s-times}{\underbrace{B_{k,n_{1}}^{\left( \alpha \right)
}\left( \xi ,q\right) B_{k,n_{2}}^{\left( \alpha \right) }\left( \xi
,q\right) ...B_{k,n_{s}}^{\left( \alpha \right) }\left( \xi ,q\right) }}%
q^{-\beta \xi }d\mu _{-q^{\beta }}\left( \xi \right) \\
&=&\dprod\limits_{i=1}^{s}\binom{n_{i}}{k}\int_{%
\mathbb{Z}
_{p}}\left[ \xi \right] _{q^{\alpha }}^{sk}\left[ 1-\xi \right] _{q^{-\alpha
}}^{n_{1}+n_{2}+...+n_{s}-sk}q^{-\beta \xi }d\mu _{-q^{\beta }}\left( \xi
\right) \\
&=&\dprod\limits_{i=1}^{s}\binom{n_{i}}{k}\sum_{l=0}^{sk}\binom{sk}{l}\left(
-1\right) ^{l+sk}\int_{%
\mathbb{Z}
_{p}}q^{-\beta \xi }\left[ 1-\xi \right] _{q^{-\alpha
}}^{n_{1}+n_{2}+...+n_{s}-sk}d\mu _{-q^{\beta }}\left( \xi \right) \\
&=&\left\{ \QATOPD. . {\left[ 2\right] _{q^{\beta }}+q^{\alpha -\beta }\frac{%
g_{n_{1}+n_{2}+...+n_{s}+1,q^{-1}}^{\left( \alpha ,\beta \right) }}{%
n_{1}+n_{2}+...+n_{s}+1},\text{ \ \ \ \ \ \ \ \ \ \ \ \ \ \ \ \ \ \ \ \ \ \
\ \ \ \ \ \ \ \ \ \ \ \ \ \ \ \ \ \ \ \ \ \ \ \ \ if }k=0,}{\dprod%
\limits_{i=1}^{s}\binom{n_{i}}{k}\sum_{l=0}^{sk}\binom{sk}{l}\left(
-1\right) ^{sk+l}\left( \left[ 2\right] _{q^{\beta }}+q^{\alpha -\beta }%
\frac{g_{n_{1}+n_{2}+...+n_{s}-l+1,q^{-1}}^{\left( \alpha ,\beta \right) }}{%
n_{1}+n_{2}+...+n_{s}-l+1}\right) ,\text{ \ \ if }k\neq 0.}\right.
\end{eqnarray*}

So from above, we have the following Theorem:

\begin{theorem}
For $s\in 
\mathbb{N}
$ with $s\geq 2$, let $n_{1},n_{2},...,n_{s},k\in 
\mathbb{N}
^{\ast }$ and $\alpha \in 
\mathbb{N}
$ with $\sum_{l=1}^{s}n_{l}>sk$. Then we have%
\begin{eqnarray*}
&&\int_{%
\mathbb{Z}
_{p}}q^{-\beta \xi }\dprod\limits_{i=1}^{s}B_{k,n_{i}}^{\left( \alpha
\right) }\left( \xi \right) d\mu _{-q}\left( \xi \right) \\
&=&\left\{ \QATOPD. . {\left[ 2\right] _{q^{\beta }}+q^{\alpha -\beta }\frac{%
g_{n_{1}+n_{2}+...+n_{s}+1,q^{-1}}^{\left( \alpha ,\beta \right) }}{%
n_{1}+n_{2}+...+n_{s}+1},\text{ \ \ \ \ \ \ \ \ \ \ \ \ \ \ \ \ \ \ \ \ \ \
\ \ \ \ \ \ \ \ \ \ \ \ \ \ \ \ \ \ \ \ \ \ if }k=0,}{\dprod\limits_{i=1}^{s}%
\binom{n_{i}}{k}\sum_{l=0}^{sk}\binom{sk}{l}\left( -1\right) ^{sk+l}\left( %
\left[ 2\right] _{q^{\beta }}+q^{\alpha -\beta }\frac{%
g_{n_{1}+n_{2}+...+n_{s}-l+1,q^{-1}}^{\left( \alpha ,\beta \right) }}{%
n_{1}+n_{2}+...+n_{s}-l+1}\right) ,\text{ \ \ if }k\neq 0.}\right.
\end{eqnarray*}
\end{theorem}

From the definition of weighted $q$-Bernstein polynomials and the binomial
theorem, we easily get%
\begin{eqnarray}
&&\int_{%
\mathbb{Z}
_{p}}\underset{s-times}{q^{-\beta \xi }\underbrace{B_{k,n_{1}}^{\left(
\alpha \right) }\left( \xi ,q\right) B_{k,n_{2}}^{\left( \alpha \right)
}\left( \xi ,q\right) ...B_{k,n_{s}}^{\left( \alpha \right) }\left( \xi
,q\right) }}d\mu _{-q^{\beta }}\left( \xi \right)  \notag \\
&=&\dprod\limits_{i=1}^{s}\binom{n_{i}}{k}\sum_{l=0}^{n_{1}+...+n_{s}-sk}%
\binom{\sum_{d=1}^{s}\left( n_{d}-k\right) }{l}\left( -1\right) ^{l}\int_{%
\mathbb{Z}
_{p}}q^{-\beta \xi }\left[ \xi \right] _{q^{\alpha }}^{sk+l}d\mu _{-q^{\beta
}}\left( \xi \right)  \notag \\
&=&\dprod\limits_{i=1}^{s}\binom{n_{i}}{k}\sum_{l=0}^{n_{1}+...+n_{s}-sk}%
\binom{\sum_{d=1}^{s}\left( n_{d}-k\right) }{l}\left( -1\right) ^{l}\frac{%
g_{l+sk+1,q}^{\left( \alpha ,\beta \right) }}{l+sk+1}.  \label{equation 17}
\end{eqnarray}

Therefore, from (\ref{equation 17}) and Theorem 8, we get interesting
Corollary as follows:

\begin{corollary}
For $s\in 
\mathbb{N}
$ with $s\geq 2$, let $n_{1},n_{2},...,n_{s},k\in 
\mathbb{N}
^{\ast }$ and $\alpha \in 
\mathbb{N}
$ with $\sum_{l=1}^{s}n_{l}>sk.$ We have 
\begin{eqnarray*}
&&\sum_{l=0}^{n_{1}+...+n_{s}-sk}\binom{\sum_{d=1}^{s}\left( n_{d}-k\right) 
}{l}\left( -1\right) ^{l}\frac{g_{l+sk+1,q}^{\left( \alpha ,\beta \right) }}{%
l+sk+1} \\
&=&\left\{ \QATOPD. . {\left[ 2\right] _{q^{\beta }}+q^{\alpha -\beta }\frac{%
g_{n_{1}+n_{2}+...+n_{s}+1,q^{-1}}^{\left( \alpha ,\beta \right) }}{%
n_{1}+n_{2}+...+n_{s}+1},\text{ \ \ \ \ \ \ \ \ \ \ \ \ \ \ \ \ \ \ \ \ \ \
\ \ \ \ \ \ \ \ if }k=0,}{\sum_{l=0}^{sk}\binom{sk}{l}\left( -1\right)
^{sk+l}\left( \left[ 2\right] _{q^{\beta }}+q^{\alpha -\beta }\frac{%
g_{n_{1}+n_{2}+...+n_{s}-l+1,q^{-1}}^{\left( \alpha ,\beta \right) }}{%
n_{1}+n_{2}+...+n_{s}-l+1}\right) ,\text{ \ \ if }k\neq 0.}\right.
\end{eqnarray*}
\end{corollary}

\end{document}